\documentclass[12pt]{article}
\usepackage{amsxtra,amssymb,amsthm,amsmath,latexsym}

\theoremstyle{plain}

\def\nd{\noindent}

\def\R{{\mathbb R}}

\def\oH{\buildrel\circ\over H}
\def\oH1{\buildrel\circ\over H\kern-.02in{}^1}
\def\qed{{\hfill $\Box$}}
\def\l{\ell}

\begin{document}

\title{ Modified Rayleigh Conjecture for static problems
   \thanks{key words: Modified Rayleigh Conjecture, 
boundary-value problems
    }
   \thanks{AMS subject classification: 35R30 }
}

\author{
A.G. Ramm\\
Mathematics Department, 
Kansas State University, \\
 Manhattan, KS 66506-2602, USA\\
ramm@math.ksu.edu\\
}

\date{}

\maketitle\thispagestyle{empty}

\begin{abstract}
Modified Rayleigh conjecture (MRC) in scattering theory was proposed and
justified by the author (J.Phys A, 35 (2002), L357-L361). MRC allows one 
to develop efficient numerical algorithms for solving 
boundary-value problems. It gives an error 
estimate for solutions. In this paper the MRC is formulated and proved for
static problems.

\end{abstract}


\section{Introduction}
Consider a bounded domain $D \subset \R^n$, $n = 3$ with a boundary $S$. The
exterior domain is $D^\prime = \R^3 \backslash D$. Assume that $S$ is
Lipschitz.
Let $ S^2$ denotes the unit sphere in $\R^3$. 
Consider the problem:
$$\nabla^2 v=0 \hbox{\ in\ } D^\prime, \quad
  v = f \hbox{\ on\ } S, \eqno{(1.1)}$$
$$v:=O(\frac 1{r})
 \quad r:=|x| \to \infty. 
\eqno{(1.2)}$$
Let $\frac x r:=\alpha\in S^2$. Denote
by $Y_\l (\alpha)$  the orthonormal spherical harmonics,
$Y_\l = Y_{\l m}, -\l \leq m \leq \l$. Let $h_\l :=\frac 
{Y_\l(\alpha)}{r^{\l 
+1}}$, $\l\geq 0,$ be  harmonic functions in $D'$.
Let the ball $B_R := \{x : |x| \leq R\}$ contain $D$.

In the region $r> R$ the solution to (1.1) - (1.2) is:
$$v(x) = 
\sum^\infty_{\l =0} c_\l h_\l, \quad
\quad r > R,\quad  
\eqno{(1.3)}$$
the summation in (1.3) and below includes summation with respect to $m$, 
$-\l \leq m \leq \l$,
and $c_\l $ are some coefficients determined by $f$.

The series (1.3)  in general does not converge up to the 
boundary $S$.
Our aim is to give a formulation of an analog of the Modified Rayleigh 
Conjecture (MRC) from [1],
which can be used in numerical solution of the 
boundary-value  problems. 
The author hopes that the MRC method
for static problems can be used as a basis for an efficient
numerical algorithm for solving boundary-value problems for Laplace
equations in domains with complicated boundaries. In [4] such
an algorithm was developed on the basis of MRC  for solving
boundary-value problems for the Helmholtz equation.
Although the boundary integral equation methods and finite elements 
methods are widely and successfully used for solving these problems, the 
method, based on 
MRC, proved to be competitive and often superior to the currently 
used methods. 

We discuss the Dirichlet condition but 
a similar argument is applicable to the Neumann and Robin boundary 
conditions. Boundary-value problems and scattering problems in rough 
domains were studied in [3]. 
  
Let us present the basic results on which the MRC method is based.

Fix $\epsilon >0$, an arbitrary small number.

{\bf Lemma 1.1.} {\it There exist $L=L(\epsilon)$ and 
$c_\l=c_\l(\epsilon)$
such that }
$$ ||\sum_{\l=0}^{L(\epsilon)}c_\l(\epsilon)h_\l -f||_{L^2(S)} \leq 
\epsilon.
\eqno{(1.4)}$$ 

If (1.4) and the boundary condition (1.1) hold, then
$$ ||v_{\epsilon}-v||_{L^2(S)}\leq \epsilon,  \quad 
v_{\epsilon}:=\sum_{\l=0}^{L(\epsilon)}c_\l(\epsilon)h_\l.
\eqno{(1.5)}$$

{\bf Lemma 1.2.} {\it If (1.4) holds then
$$ ||v_{\epsilon}-v||=O(\epsilon) \quad \epsilon \to 0, \quad
\eqno{(1.6)}$$
where $||\cdot||:= ||\cdot||_{H_{loc}^m(D')}+||\cdot||_{L^2(D'; 
(1+|x|)^{-\gamma})}$, $\gamma >1$, $m>0$ is an arbitrary integer,
and $H^m$ is the Sobolev space.} 

In particular, (1.6) implies
$$ ||v_{\epsilon}-v||_{L^2(S_R)}=O(\epsilon) \quad  \epsilon \to 0. 
\eqno{(1.7)}$$

Let us formulate an analog of the Modified Rayleigh Conjecture (MRC):  

{\bf Theorem 1 (MRC):} {\it For an arbitrary small $\epsilon>0$ there 
exist
$L(\epsilon)$ and $c_\l(\epsilon), 0\leq \l \leq L(\epsilon)$,
such that (1.4) and (1.6)  hold.}

Theorem 1 follows from  Lemmas 1.1 and 1.2.

For the Neumann boundary condition one minimizes
$ ||\frac {\partial [\sum_{\l=0}^{L}c_\l\psi_\l]}{\partial 
N}-f||_{L^2(S)}$
with respect to $c_\l$. Analogs of Lemmas 1.1-1.2 are valid and their 
proofs are essentially the same.

If the boundary data $f\in C(S)$, then one can use $C(S)-$ norm 
in (1.4)-(1.7), and an analog of Theorem 1 then 
follows immediately from the maximum principle.

In Section 2 we discuss the usage of MRC in solving boundary-value 
problems. In Section 3 proofs are given.

\section{Solving boundary-value problems by MRC.} 
To solve problem (1.1)-(1.2)
using MRC, fix a small $\epsilon >0$
and find $L(\epsilon)$ and $c_\l(\epsilon)$ such that 
(1.4) holds. This is possible by Lemma 1.1 and can be done numerically
by minimizing $||\sum_0^Lc_\l h_\l -f||_{L^2(S)}:=\phi 
(c_1,.....,c_L)$. If the minimum of $\phi$ is larger than $\epsilon$, 
then increase $L$ and repeat the minimization. Lemma 1.1 guarantees the
existence of such $L$ and $c_\l$ that the minimum is less than $\epsilon$.
Choose the smallest $L$ for which this happens and
define $v_\epsilon:=\sum^L_{\l = 0} c_\l h_\l$. Then, by Lemma 1.2,
$v_\epsilon$ is the approximate solution to problem (1.1)-(1.2)
with the accuracy $O(\epsilon)$ in the norm $||\cdot||$. 

\section{Proofs.}

{\bf Proof of Lemma 1.1.}   We start with the claim:

{\bf Claim}: 
  {\it the restrictions of harmonic functions $h_\l$ on $S$
form a total set in $L^2(S)$.} 

Lemma 1.1 follows from this claim.
Let us prove the claim. Assume the contrary. Then there is
a function $g\neq 0$ such that $\int_S g(s)h_\l(s) ds=0\,\, \forall \l\geq 
0.$ This implies $V(x):=\int_S g(s)|x-s|^{-1}ds=0\,\, \forall x\in D'$. 
Thus $V=0$ on $S$, and since $\Delta V=0$ in $D$, one concludes that $V=0$
in $D$. Thus $g=0$ by the jump formula for the normal derivatives of
the simple layer potential $V$. This contradiction  proves
the claim. Lemma 1.1 is proved. 
$\Box$

\nd{\bf Proof of Lemma 1.2.} 
By Green's formula one has
$$ w_\epsilon (x) = \int_S w_\epsilon (s) G_N (x,s) ds, \quad
  \| w_\epsilon  \|_{L^2(S)} < \epsilon,\quad  w_\epsilon:=v_\epsilon-v. 
\eqno{(3.1)}$$
Here $N$ is the unit normal to $S$, pointing into $D'$, and $G$ is the 
Dirichlet Green's function of the Laplacian in $D^\prime$:
$$\nabla^2 G =-\delta (x-y) \hbox{\ in\ } D^\prime,
  \quad G=0 \hbox{\ on\ } S, \eqno{(3.2)}$$
$$ G = O(\frac 1 r), \quad r\to \infty. \eqno{(3.3)}$$
From (3.1) one gets (1.7) and (1.6) with 
$H^m_{loc}(D')-$norm immediately by the 
Cauchy inequality. Estimate (1.6) in the region 
$B_{R}^\prime:=\R^3\setminus B_R$ follows
from the estimate
$$\left|G_N (x,s)\right| \leq \frac{c}{1 + |x|}, \quad |x| 
\geq R.
  \eqno{(3.4)}$$
In the region $B_R\backslash D$ estimate (1.6) follows from local elliptic 
estimates for $w_\epsilon:=v_\epsilon -v$, which 
imply that
$$\|w_\epsilon \|_{L^2(B_R\backslash D)} \leq c\epsilon. 
\eqno{(3.5)}$$
Let us recall the elliptic estimate we have used.
Let  $D'_{R}:=B_R\backslash D$ and $S_R$ be the boundary of $B_R$.
 Recall the elliptic estimate  for the solution to
homogeneous Laplace equation in $D'_{R}$ ( see \cite{2}, p.189):
$$\|w_\epsilon \|_{H^{0.5}(D'_{R})} \leq c [ 
||w_\epsilon||_{L^2(S_R)} + ||w_\epsilon||_{L^2(S)}].
\eqno{(3.6)}$$ 
The estimates 
$||w_\epsilon||_{L^2(S_R)}=O(\epsilon)$, 
 $||w_\epsilon||_{L^2(S)}=O(\epsilon)$, and (3.6) yield
(1.6). Lemma 1.2 is proved. \qed

\end{document}